\documentclass[11pt]{article}
\usepackage{amsmath,amssymb,latexsym}

\newtheorem{thm}{Theorem}
\newtheorem{lem}[thm]{Lemma}
\newtheorem{cor}[thm]{Corollary}
\newtheorem{rmk}[thm]{Remark}

\newtheorem{exm}{Example}
\newcommand{\bea}{\begin{eqnarray}}
\newcommand{\eea}{\end{eqnarray}}
\newcommand{\bean}{\begin{eqnarray*}}
\newcommand{\eean}{\end{eqnarray*}}
\newcommand{\beq}{\begin{equation}}
\newcommand{\eeq}{\end{equation}}
\newcommand{\bac}{\begin{array}{c}}
\newcommand{\ball}{\begin{array}{ll}}
\newcommand{\ea}{\end{array}}

\def\({\left(}
\def\){\right)}

\title{On Smooth Whitney Extensions of almost isometries with small distortion, Interpolation and Alignment in $\Bbb R^D$-Part 1}
\author{Steven B. Damelin \thanks{Zentralblatt MATH, FIZ Karlsruhe – Leibniz Institute\, email: steve.damelin@gmail.com}\and Charles Fefferman, \thanks{Department of Mathematics; Fine Hall, Washington Road, Princeton NJ 08544-1000 USA,\, email: cf@math.princeton.edu} }
\date{}

\begin{document}

\maketitle

\thispagestyle{empty}
\parskip=10pt

\begin{abstract}
In this paper,  we study the following problem: Let $D\geq 2$ and let $E\subset \mathbb R^D$ be finite satisfying certain conditions. Suppose that we are given a map $\phi:E\to \mathbb R^D$ with
$\phi$ a small distortion on $E$. How can one decide whether $\phi$ extends to a smooth small distortion $\Phi:\mathbb R^D\to \mathbb R^D$ which agrees with $\phi$ on $E$. We also ask how to decide if in addition $\Phi$ can be approximated well by certain rigid and non-rigid motions from $\mathbb R^D\to \mathbb R^D$. Since $E$ is a finite set, this question is basic to interpolation and alignment of data in $\mathbb R^D$.  The work in this paper appears in the research memoir \cite{SDam}. 
\medskip

\noindent AMS-MSC Classification: 58C25, 42B35, 94A08, 94C30, 41A05, 68Q25, 30E05, 26E10, 68Q17.

\end{abstract}

\vskip 0.2in
\tableofcontents

\section{Introduction.}
\setcounter{thm}{0}
\setcounter{equation}{0}

We will work here and throughout in Euclidean space $\mathbb R^D$, $D\geq 2$ and $|.|$ will always denote the Euclidean norm on $\mathbb R^D$.  Henceforth, $D$ is chosen and fixed. 

We study the following problem:
\medskip

(P)\, {\bf Problem.}\, Let $E\subset \mathbb R^D$ be finite satisfying certain conditions which we will define precisely in the statements of our main results. Suppose that we are given a map $\phi:E\to \mathbb R^D$ with
$\phi$ a small distortion on $E$. How can one decide whether $\phi$ extends to a smooth small distortion $\Phi:\mathbb R^D\to \mathbb R^D$ which agrees with $\phi$ on $E$. We also ask how to decide if $\Phi$ can be approximated well by certain rigid and non-rigid motions from $\mathbb R^D\to \mathbb R^D$.
\medskip

\subsection{Some insight}

We provide some brief insight into problem (P):

\begin{itemize}
\item  {\bf Example Whitney problem 1.} 
\item Suppose we are given two sets of distinct points in $\mathbb R^D$. The data sets are indexed by the same set (labeled):  
\item Pairwise distances between corresponding points are equal in the two data sets.
\item  In other words, the sets are isometric.
\item  Can this correspondence be extended to an isometry of the whole of $\mathbb R^D$?
\end{itemize}

\begin{itemize}
\item {\bf Example Whitney extension problem 2.}
\item $E\subset \mathbb R^D$ be finite with some geometry.  
\item $\phi:E\to \mathbb R^D$, $\phi$ a small distortion on $E$. How can one decide whether $\phi$ extends to a small distorted diffeomorphism (in particular 1-1 and onto)   $\Phi:\mathbb R^D\to \mathbb R^D$ which agrees with $\phi$ on $E$. 
\item How to decide if $\Phi$ can be approximated well by certain rigid and non-rigid motions from $\mathbb R^D\to \mathbb R^D$.  

\end{itemize}

\begin{itemize}
\item {\bf Example Whitney extension problem 3.}
\item Let $1\leq m\leq \infty$, $U\subset\mathbb R^D$ open.
\item Let $\phi:U\to \mathbb R^D$ be a $C^m$ map. If $E\subset U$ is a "fat" compact set  (with some geometry) and the restriction of $\phi$ to $E$ is an almost isometry 
with small distortion,  how to decide when there exists a $C^m(\mathbb R^D)$ one-to-one and onto almost isometry $\Phi:\mathbb R^D\to \mathbb R^D$ with smooth distortion
which agrees with $\phi$ in a neighborhood of $E$ and a Euclidean motion $A:\mathbb R^D\to \mathbb R^D$ away from $E$.  
\end{itemize}

Given $E$ is finite, Problems 1 and 2 are interpolation problems.

\begin{itemize}
\item {\bf Whitney extension problem 1}: The answer to this problem has been well known for sometime, see for example \cite{WW1}.
\item Let $y_1,...y_k$ and $z_1,...z_k$ be distinct points in $\mathbb R^D$.  Suppose that
\[
|z_i-z_j|=|y_i-y_j|,\,\,  1\leq i,j\leq k,\, i\neq j.
\]
Then there exists a smooth isometry $\Psi:\mathbb R^D\to \mathbb R^D$, i.e
\[
\frac{|\Psi(x)-\Psi(y)|}{|x-y|}=1,\, \, x,y\in \mathbb R^D,\, x\neq y
\]
satisfying
\[
\Psi(y_i)=z_i,\, \, 1\leq i\leq k.
\]
\end{itemize}

\begin{itemize}
\item {\bf Example Whitney extension problem 4.}
\item $E\subset \mathbb R^D$ arbitrary. 
\item $\phi:E\to \mathbb R$.  How can one decide whether $\phi$ extends to a smooth map $\Phi:\mathbb R^D\to \mathbb R$ which agrees with $\phi$ on $E$. 
\item For $D=1$ and $E$ compact, this is the well-known Whitney extension problem. Progress on this problem has been made by G. Glaeser \cite{Gl}, Y. Brudnyi and P. Shvartsman and his collaborators, see  \cite{B,BS1,BS2,BS3,BS4} and in addition many papers (P. Shvartsman, arxiv), E. Bierstone, P. Milman and W. Pawlucki \cite{BMP} and  C. Wells \cite{WW1}. See also Zobin \cite{Z,Z1}, E. Le Gruyer \cite{LG} and  C. Fefferman and his collaborators,  see \cite{F2,F3,F5} and in addition, many papers  (C. Fefferman, arxiv).
\end{itemize}

\begin{itemize}
\item {\bf  Shape registration or allignment problem}.
\item Given two distinct configurations in $\Bbb R^D$, do there exist combinations of rotations, translations, reflections and compositions of these which map the one configuration onto the other. This is called the shape registration or allignment problem where. In this sense, an old result called the 
Procrustes problem see for example \cite{WW}, establishes:
\item Let $y_1,...y_k$ and $z_1,...z_k$ be two sets of distinct points in $\mathbb R^D$. Suppose that
\[
|z_i-z_j|=|y_i-y_j|,\,  1\leq i,j\leq k,\, i\neq j.
\]
Then there exists a Euclidean motion motion $\Phi_0:x\to Tx+x_0$ such that $\Phi_0(y_i)=z_i,\, i=1,...,k.$. If $k\leq D$, then $\Phi_0$ is proper.
\end{itemize}

\section{Main results}

\subsection{Notation}

Now and throughout, given a positive number $\varepsilon$, we assume that $\varepsilon$ is less than or equal to a positive constant determined by $D$. Henceforth, $C,c,c_1,...$ will denote positive controlled constants depending only on $D$. These symbols need not denote the same constant in different occurrences. Throughout, we will use other symbols to denote positive constants which may not depend on $D$ or depend on $D$ but also on other quantities. These symbols need not denote the same constant in different occurrences.

Following are our main results:

\subsection{Main Results}

\begin{thm}
Given $\varepsilon>0$, there exists $\delta>0$ depending on $\varepsilon$ such that the following holds: Let $y_1,...y_k$ and $z_1,...z_k$ be two $1\leq k\leq D$ sets of distinct points in $\mathbb R^D$. Suppose that
\beq
(1+\delta)^{-1}\leq \frac{|z_i-z_j|}{|y_i-y_j|}\leq (1+\delta),\, 1\leq i,j\leq k,\, i\neq j.
\label{e:emotionsa}
\eeq
Then there exists a diffeomorphism (in particular 1-1 and onto) $\Phi:\mathbb R^D\to \mathbb R^D$ with
\beq
(1+\varepsilon)^{-1}|x-y|\leq |\Phi(x)-\Phi(y)|\leq |x-y|(1+\varepsilon),\, x,y\in \mathbb R^D
\label{e:emotionsaa}
\eeq
satisfying
\beq
\Phi(y_i)=z_i,\, 1\leq i\leq k.
\label{e:emotionsaaa}
\eeq
\label{t:lemmafive}
\end{thm}

\begin{thm}
Let $\varepsilon>0$ and let $y_1,...y_k$ and $z_1,...z_k$ be two $1\leq k\leq D$ sets of distinct points in $\mathbb R^D$ with
\beq
\sum_{i\neq j}|y_i-y_j|^2+\sum_{i\neq j}|z_i-z_j|^2=1,\, y_1=z_1=0.
\label{e:lemmafourpoints}
\eeq
Then there exist constants $c,c_1'>0$ depending only on $D$ such the following holds: Set $\delta=c_1\varepsilon'^{c}$ and suppose 
\beq
||z_i-z_j|-|y_i-y_j||<\delta, \, i\neq j.
\label{e:lemmafourpointsa}
\eeq
Then there exists a diffeomorphism $\Phi:\mathbb R^n\to \mathbb R^D$ satisfying $(\ref{e:emotionsaa})$ and $(\ref{e:emotionsaaa})$.
\label{t:lemmafiveprime}
\end{thm}

\begin{thm}
Given $\varepsilon>0$, there exists $\delta>0$ depending on $\varepsilon$, such that the following holds. Let $y_1,...,y_k$ and
$z_1,...,z_k$ be two $k\geq 1$ point configurations in $\mathbb R^D$. Suppose that 
\[
(1+\delta)^{-1}\leq \frac{|z_i-z_j|}{|y_i-y_j|}\leq 1+\delta,\, i\neq j.
\]
Then, there exists a Euclidean motion $\Phi_0:x\to Tx+x_0$ such that
\beq
|z_i-\Phi_0(y_i)|\leq \varepsilon {\rm diam}\left\{y_1,...,y_k\right\}
\label{e:emotionsb}
\eeq
for each $1\leq i\leq k$. If $k\leq D$, then we can take $\Phi_0$ to be a proper Euclidean motion on $\mathbb R^{D}$.
\label{t:lemmafour}
\end{thm}
\smallskip

\begin{thm}
Let $\varepsilon>0$ and let $y_1,...y_k$ and $z_1,...z_k$ be two $k\geq 1$ sets of distinct points in $\mathbb R^D$ with
\[
\sum_{i\neq j}|y_i-y_j|^2+\sum_{i\neq j}|z_i-z_j|^2=1,\, y_1=z_1=0.
\]
Then there exist constants $c,c_1'>0$ depending only on $D$ such the following holds: Set $\delta=c_1\varepsilon'^{c}$ and suppose 
\[
||z_i-z_j|-|y_i-y_j||<\delta, \, i\neq j.
\]
Then, there exists a Euclidean motion $\Phi_0:x\to Tx+x_0$ such that
\beq
|z_i-\Phi_0(y_i)|\leq \varepsilon
\label{e:emotionsu}
\eeq
for each $1\leq i\leq k$. If $k\leq D$, then we can take $\Phi_0$ to be a proper Euclidean motion on $\mathbb R^{D}$.
\label{t:lemmafoura}
\end{thm}

The reminder of this paper proves Theorem~\ref{t:lemmafour}, Theorem~\ref{t:lemmafoura}, Theorem~\ref{t:lemmafive} and Theorem~\ref{t:lemmafiveprime}.  We begin with:

\section{$\varepsilon$-distorted diffeomorphisms in $\mathbb R^D$.}
\setcounter{equation}{0}

Let $\varepsilon>0$ and $\Phi:\mathbb R^D\to \mathbb R^D$ be a diffeomorphism. We say that $\Phi$ is ``$\varepsilon$-distorted" provided
\[
(1+\varepsilon)^{-1}I\leq (\Phi'(x))^{T}\Phi'(x)\leq (1+\varepsilon)I
\]
as matrices, for all $x\in \mathbb R^D$. Here, $I$ denotes the identity matrix in $\mathbb R^D$.
\medskip

We will use the following properties of $\varepsilon$-distorted maps:
\begin{itemize}
\item If $\Phi$ is $\varepsilon$-distorted and $\varepsilon<\varepsilon'$, then $\Phi$ is $\varepsilon'$-distorted.
\item If $\Phi$ is $\varepsilon$-distorted, then so is $\Phi^{-1}$
\item If $\Phi$ and $\Psi$ are $\varepsilon$-distorted, then the composition map $\Phi o\Psi$ is $C\varepsilon$-distorted.
\item Suppose $\Phi$ is $\varepsilon$-distorted. If $\tau$ is a piecewise smooth curve in $\mathbb R^{D}$, then the length of
$\Phi(\tau)$ differs from that of $\tau$ by at most a factor of $(1+\varepsilon)$. Consequently, if $x,x'\in \mathbb R^D$, then $|x-x'|$ and $|\Phi(x)-\Phi(x')|$ differ by at most a factor $(1+\varepsilon)$.
\end{itemize}

The first three properties follow easily from the definition. The fourth property follows for example using Bochner's theorem and it is this property which implies immediately that every $\varepsilon$-distorted diffeomorphism $\Phi:\mathbb R^D\to \mathbb R^D$ satisfies,
\beq
(1+\varepsilon)^{-1}|x-y|\leq |\Phi(x)-\Phi(y)|\leq |x-y|(1+\varepsilon),\, x,y\in \mathbb R^D.
\label{e:edistortion}
\eeq 
It is also not difficult to show that (\ref{e:edistortion}) together with the fact that $(\Phi'(x))^{T}\Phi'(x)$ is real and symmetric implies that 
\[
|(\Phi'(x))^{T}\Phi'(x)-I|\leq C\varepsilon,\, x\in \mathbb R^D.
\]
This follows from working in an orthonormal basis for which $(\Phi'(x))^{T}\Phi'(x)$ is diagonal.
\subsection{Examples.}
\setcounter{equation}{0}

In this subsection, we provide two useful examples of $\varepsilon$-distorted diffeomorphisms of $\mathbb R^D$ which we will use.

\subsection{Slow twists.}

We are interested in rotations which are $\varepsilon$-distorted diffeomorphisms, rotate non-rigidly close to the origin and become more rigid away. See Lemma~\ref{l:lemmaone}.

\begin{exm}
{\rm Let $\varepsilon>0$ and $x\in \mathbb R^D$. Let $S(x)$ be the $D\times D$ block-diagonal matrix
\[
\left(
\begin{array}{llllll}
D_1(x) & 0 & 0 & 0 & 0 & 0 \\
0 & D_2(x) & 0 & 0 & 0 & 0 \\
0 & 0 & . & 0 & 0 & 0 \\
0 & 0 & 0 & . & 0 & 0 \\
0 & 0 & 0 & 0 & . & 0 \\
0 & 0 & 0 & 0 & 0 & D_r(x)
\end{array}
\right)
\]
where for each $i$ either $D_i(x)$ is the $1\times 1$ identity matrix or else
\[
D_i(x)=\left(
\begin{array}{ll}
\cos f_i(|x|) & \sin f_i(|x|) \\
-\sin f_i(|x|) & \cos f_i(|x|)
\end{array}
\right)
\]
where $f_i:\mathbb R\to \mathbb R$ are functions (corresponding to the blocks $D_i$) satisfying the condition: $t|f_i'(t)|<c\varepsilon$ some $c>0$ small enough, uniformly for $t\geq 0$. Let 
$\Phi(x)=\Theta^{T}S(\Theta x)$ where $\Theta$
is any fixed matrix in $SO(D)$. Then it is straightforward to check that $\Phi$ is a $\varepsilon$-distorted diffeomorphism and we call it a {\it slow twist} (in analogy to rotations). }
\label{e:Example1}
\end{exm}
. 

\subsection{Slides.}

We are interested in translations which are $\varepsilon$-distorted diffeomorphisms and translate non-rigidly close to the origin and become more rigid away. See Lemma~\ref{l:lemmatwo}.

\begin{exm}
{\rm Let $\varepsilon>0$ and let $g:\mathbb R^{D}\to \mathbb R^{D}$
be a smooth map such that $|g'(t)|<c\varepsilon$ some $c>0$ small enough, uniformly for $t\in \mathbb R^D$. 
Consider the map $\Phi(t)=t+g(t)$, $t\in \mathbb R^D$.   Then $\Phi$ is a $\varepsilon$-distorted diffeormorphism . We call the map $\Phi$ a {\it slide} (in analogy to translations). }
\label{e:Example2}
\end{exm}

\begin{rmk}
{\rm Lemma~\ref{l:lemmaone} and Lemma~\ref{l:lemmatwo} below show us how to approximate $SO(D)$ and proper Euclidean motions by $\varepsilon$-distorted diffeomorphisms. Conversely, it is not too difficult to show, see \cite{DF4} 
that given $\varepsilon>0$ and $\Phi:\mathbb R^D\to \mathbb R^D$ a $\varepsilon>0$-distorted diffeomorphism, there exists a Euclidean motion $\Phi_0:\mathbb R^D\to \mathbb R^D$ with $|\Phi(x)-\Phi_0(x)|\leq C\varepsilon$ for $x\in B(0,10)$ the ball in $\mathbb R^D$ with center 0 and radius 10.
Actually, using the well-known John-Nirenberg inequality, in \cite{DF4} we proved the following:}
\medskip

{\bf BMO theorem for $\varepsilon$-distorted diffeomorphisms}:\, Let $\varepsilon>0$, $\Phi:\mathbb R^D\to \mathbb R^D$, a $\varepsilon$-distorted diffeomorphism and let $B\in \mathbb R^D$ be a ball. There exists $T_B\in O(D)$ and 
$C>0$ such that for all $\lambda\geq 1$, we have
\[
{\rm vol}\left\{x\in B:|\Phi'(x)-T_B(x)|>C\varepsilon\lambda\right\}\leq \exp(-\lambda){\rm vol}(B)
\]
{\rm and slow twists in Example~\ref{e:Example1} show that the estimate above is sharp. The set 
\[
\left\{x\in B:|\Phi'(x)-T_B(x)|\geq C\varepsilon\lambda\right\}
\] may well be small in a more refined sense but we did not pursue this investigation in \cite{DF4}}.
\end{rmk}

\section{Approximation by Euclidean motions in $\mathbb R^D$.}
\setcounter{equation}{0}

We begin with the proof of Theorem~\ref{t:lemmafour}.
\medskip

\noindent

{\bf Proof} \ Suppose not. Then for each $l\geq 1$, we can find points $y_1^{(l)},...,y_k^{(l)}$ and $z_1^{(l)},...,z_k^{(l)}$ in
$\mathbb R^D$ satisfying (\ref{e:emotionsa}) with $\delta=1/l$ but not satisfying (\ref{e:emotionsb}). Without loss of generality, we may suppose that ${\rm diam}\left\{y_1^{(l)},...,y_k^{(l)}\right\}=1$ for each $l$ and that $y_1^{(l)}=0$ and
$z_l^{(1)}=0$ for each $l$. Thus $|y_i^{(l)}|\leq 1$ for all $i$ and $l$ and
\[
(1+1/l)^{-1}\leq \frac{|z_i^{(l)}-z_j^{(l)}|}{|y_i^{(l)}-y_j^{(l)}|}\leq (1+1/l)
\]
for $i\neq j$ and any $l$.
However, for each $l$, there does not exist an Euclidean motion
$\Phi_0$ such that
\beq
|z_i^{(l)}-\Phi_0(y_i^{(l)})|\leq \varepsilon
\label{e:emotionsc}
\eeq
for each $i$. Passing to a subsequence, $l_1,l_2,l_3,...,$ we may assume
\[
y_i^{(l_{\mu})}\to y_i^{\infty},\, \mu\to \infty
\]
and
\[
z_i^{(l_{\mu})}\to z_i^{\infty},\, \mu\to \infty.
\]
Here, the points $y_i^{\infty}$ and $z_i^{\infty}$ satisfy
\[
|z_i^{\infty}-z_j^{\infty}|=|y_i^{\infty}-y_j^{\infty}|
\]
for $i\neq j$. Hence, by Theorem~\ref{t:Theorem 1}, there is an Euclidean motion $\Phi_0:\mathbb R^D\to \mathbb R^D$ such that $\Phi_0(y_i^{\infty})=z_i^{\infty}$. Consequently,
for $\mu$ large enough, (\ref{e:emotionsc}) holds with $l_{\mu}$. This contradicts the fact that for each $l$, there does not exist a $\Phi_0$ satisfying (\ref{e:emotionsc}) with $l$.
Thus, we have proved all the assertions of the theorem except that we can take $\Phi_0$ to be proper if $k\leq D$. To see this, suppose that $k\leq D$ and let $\Phi_0$ be an improper Euclidean motion such that
\[
|z_i-\Phi_0(y_i)|\leq \varepsilon{\rm diam}\left\{y_1,...,y_k\right\}
\]
for each $i$. Then, there exists an improper Euclidean motion $\Psi_0$ that fixes $y_1,...,y_k$ and in place of $\Phi_0$, we may use $\Psi_0 o\Phi_0$ in the conclusion of the Theorem. The proof of the Theorem is complete. $\Box$.
\medskip

We now prove Theorem~\ref{t:lemmafoura}.
\medskip

\noindent

{\bf Proof} \ The idea of the proof relies on the following Lojasiewicz inequality, see \cite{SJS}.
Let $f:U\to \Bbb R$ be a real analytic function on an open set $U$ in $\Bbb R^D$ and $Z$ be the zero set of $f$. Assume that $Z$ is not empty. Then for a compact set $K$ in $U$, there exist positive constants $\alpha$ and $\alpha'$ depending on $f$ and $K$ such that for all $x\in K$, $|x-Z|^{\alpha}\leq \alpha'|f(x)|$.
It is easy to see that using this, one may construct approximating points $y_1',...,y_k',z_1',...,z_k'\in \mathbb R^D$ (zeroes of a suitable $f$) with the following two properties:
(1) There exist positive constants $c,c_1>0$ such that
\[
|y_i-y_i'|\leq c\varepsilon^{c_1}
\]
and
\[
|z_i-z_i'|\leq c\varepsilon^{c_1}.
\]
(2) $|y_i'-y_j'|=|z_i'-z_j'|$ for every $i,j$.
Thanks to (2), we may choose a Euclidean motion $\Phi_0:\mathbb R^D\to \mathbb R^D$ so that $\Phi_0(y^\prime_i)= z^\prime_i$ for each $i$.
Also, thanks to (1), 
\[
|\Phi_0(y_i)-\Phi_0(y_i')|\leq c_2\varepsilon^{c_3}
\]
So it follows that
\[
|\Phi_0(y_i)-z_i|\leq c_4\varepsilon^{c_5}
\]
which is what we needed to prove. $\Box$

\section{Auxiliary results.}
\setcounter{equation}{0}

We need two lemmas and a corollary which are consequences of the definitions of slow twists and slides.
\medskip

As a consequence of slow twists (see Example~\ref{e:Example1}), we have the following:

\begin{lem}
Given $\varepsilon>0$, there exists $\eta>0$ depending on $\varepsilon$ for which the following holds. Let $\Theta\in SO(D)$ and let $0<r_1\leq \eta r_2$.
Then, there exists an $\varepsilon$-distorted diffeomorphism $\Phi:\mathbb R^{D}\to \mathbb R^{D}$ such that
\[
\left\{
\begin{array}{ll}
\Phi(x)=\Theta x, & |x|\leq r_1 \\
\Phi(x)=x, & |x|\geq r_2
\end{array}
\right.
\]
\label{l:lemmaone}
\end{lem}

As a consequence of slides (see Example~\ref{e:Example2}), we have the following:

\begin{lem}
Given $\varepsilon>0$, there exists $\eta>0$ depending on $\varepsilon$ such that the following holds. Let $\Phi_1: x\to Tx+x_0:\mathbb R^D\to \mathbb R^D$ be a proper Euclidean
motion. Let $r_1,r_2>0$.
Suppose $0<r_1\leq \eta r_2$ and $|x_0|\leq c\varepsilon r_1$.
Then there exists an $\varepsilon$-distorted diffeomorphism $\Phi_2:\mathbb R^{D}\to \mathbb R^{D}$ such that
\[
\left \{
\begin{array}{ll}
\Phi_{2}(x)=\Phi_{1}(x), & |x|\leq r_1 \\
\Phi_{2}(x)=x, & |x|\geq r_2
\end{array}
\right.
\]
\label{l:lemmatwo}
\end{lem}

As a result of Lemma~\ref{l:lemmaone} and Lemma~\ref{l:lemmatwo} we have the following corollary.

\begin{cor}
Given $\varepsilon>0$, there exists $\eta>0$ depending on $\varepsilon$ such that the following holds. Let $0<r_1\leq \eta r_2$ and let $x,x'\in \mathbb R^{D}$ with $|x-x'|\leq c\varepsilon r_1$ for some $c>0$ and $|x|\leq r_1$. Then, there exists an $\varepsilon$-distorted diffeomorphism $\Phi:\mathbb R^{D}\to\mathbb R^{D}$ such that $\Phi(x)=x'$ and $\Phi(y)=y$ for $|y|\geq r_2$.
\label{c:corollaryone}
\end{cor}

We next need to introduce a technique which involves the combinatorics of certain hierarchical clusterings of finite subsets of $\mathbb R^D$. In this regard, we need the following:

\begin{lem}
Let $k\geq 2$ be a positive integer and let $0<\eta\leq 1/10$. Let $E\subset \mathbb R^D$ be a set consisting of $k$
distinct points. Then, we can partition $E$ into sets $E_1,E_2,...E_{{\mu}_{\rm max}}$ and we can find a positive integer $l$
$(10\leq l\leq 100+\binom{k}{2})$ such that the following hold:
\beq
{\rm diam}(E_{\mu})\leq \eta^{l}{\rm diam}(E)
\eeq
for each $\mu$ and
\beq
{\rm dist}(E_{\mu},E_{\mu'})\geq \eta^{l-1}{\rm diam}(E),\, {\rm for}\, \mu\neq \mu'.
\eeq
\label{l:lemmathree}
\end{lem}

\noindent
{\bf Proof} \ We define and equivalence relation on $E$ as follows. Define a relation $\sim$ on $E$ by saying that $x\sim x'$, for $x,x'\in E$ if and only if $|x-x'|\leq \eta^{l}{\rm diam}(E)$ for a fixed positive integer $l$ to be defined in a moment. By the pigeonhole principle, we can always find a positive integer $l$ such that
\[
|x-x'| \notin (\eta^l{\rm diam}(E),\eta^{l-1}{\rm diam}(E)], \, x,x'\in E.
\]
and such that $10\leq l\leq 100+\binom{k}{2}$. Let us choose and fix such an $l$ and use it for $\sim$ as defined above.
Then $\sim$ is an equivalence relation and the equivalence classes of $\sim$ partition $E$
into the sets $E_1,...,E_{{\mu}_{\rm max}}$ with the properties as required.
$\Box$

\section{A special case of Theorem~\ref{t:lemmafive}.}
\setcounter{equation}{0}

In this section, we prove a special case of Theorem~\ref{t:lemmafive}. This is given in the following theorem.

\begin{thm}
Let $\varepsilon>0$ and let $m$ be a positive integer. Let $\lambda>0$ be less than a small enough constant depending only on $\varepsilon$, $m$ and $D$. Let $\delta>0$ be less than a small enough constant depending only on
$\lambda$, $\varepsilon$, $m$ and $D$. Then the following holds: Let $E:=y_1,...y_k$ and $E':=z_1,...z_k$ be $k\geq 1$ distinct points in $\mathbb R^D$ with $k\leq D$ and $y_1=z_1$. Assume the following:
\beq
|y_i-y_j|\geq \lambda^m{\rm diam}\left\{y_1,...,y_k\right\},\, i\neq j
\label{e:exactone}
\eeq
and
\beq
(1+\delta)^{-1}\leq \frac{|z_i-z_j|}{|y_i-y_j|}\leq (1+\delta),\, i\neq j.
\label{e:exacttwo}
\eeq
Then, there exists an $\varepsilon$-distorted diffeomorphism $\Phi:\mathbb R^D\to \mathbb R^D$ such that
\begin{equation}
\Phi(y_i)=z_i, 1\leq i\leq k
\label{e:exactthree}
\eeq
and
\beq
\Phi(x)=x \ \mbox{for} \ |x-y_1|\geq \lambda^{-1/2}{\rm diam}\left\{y_1,...,y_k\right\}.
\label{e:exactfour}
\eeq
\label{t:theorem3exact}
\end{thm}

\noindent
{\bf Proof} \ Without loss of generality, we may take $y_1=z_1=0$ and ${\rm diam}\left\{y_1,...,y_k\right\}=1$.
Applying Theorem~\ref{t:lemmafour} with $10^{-9}\varepsilon \lambda^{m+5}$ in place of $\varepsilon$, we obtain a proper Euclidean motion
\beq
\Phi_0:x\to Tx+x_0
\label{e:exactfive}
\eeq
such that
\beq
|\Phi_0(y_i)-z_i|\leq 10^{-9}\varepsilon \lambda^{m+5}
\label{e:exactsix}
\eeq
for each $i$. In particular, taking $i=1$ and recalling that $y_1=z_1=0$, we find that
\beq
|x_0|\leq 10^{-9}\varepsilon \lambda^{m+5}.
\label{e:exactseven}
\eeq

For each $i$, we consider the balls
\beq
B_i=B(\Phi_0(y_i),\lambda^{m+3}),\, B_i^{+}=B(\Phi_0(y_i),\lambda^{m+1}).
\label{e:exacteight}
\eeq

Note that (\ref{e:exactone}) shows that the balls $B_i^{+}$
have pairwise disjoint closures since $\Phi_0$
is an Euclidean motion. Applying Corollary~\ref{c:corollaryone}
to Lemma~\ref{l:lemmatwo}, we obtain for each $i$,
a $\varepsilon$-distorted diffeomorphism
$\Psi_i:\mathbb R^D\to \mathbb R^D$ such that
\beq
(\Psi_i o\Phi_0)(y_i)=z_i
\label{e:exactnine}
\eeq
and
\beq
\Psi_i(x)=x
\label{e:exactten}
\eeq
outside $B_i^{+}$.
In particular, we see that
\beq
\Psi_i:B_i^+\to B_i^+
\label{e:exacteleven}
\eeq
is one to one and onto. We may patch the $\Psi_i$ together into a single map $\Psi:\mathbb R^D\to \mathbb R^D$ by setting
\beq
\Psi(x):=\left\{
\begin{array}{ll}
\Psi_i(x), & x\in B_i^{+} \\
x, & x\notin \cup_j B_j^{+}
\end{array}
\right\}.
\label{e:exacttwelve}
\eeq
Since the $B_i^{+}$ have pairwise disjoint closures (\ref{e:exactten}) and (\ref{e:exacteleven}) show that
$\Psi$ maps $\mathbb R^D$ to $\mathbb R^D$ and is one to one and onto. Moreover, since each
$\Psi_i$ is $\varepsilon$-distorted, it now follows easily that
\beq
\Psi:\mathbb R^D\to \mathbb R^D
\label{e:exactthirteen}
\eeq
is an $\varepsilon$-distorted diffeomorphism.
From (\ref{e:exacteight}), (\ref{e:exactnine}) and (\ref{e:exacttwelve}),
we also see that
\beq
(\Psi o\Phi_0)(y_i)=z_i, \forall i.
\label{e:exactforteen}
\eeq
Suppose $x\in \mathbb R^D$ with $|x|\geq 5$.
Then (\ref{e:exactfive}) and (\ref{e:exactseven}) show that
$|\Phi_0(x)|\geq 4$. On the other hand, each $y_i$ satisfies
\[
|y_i|=|y_i-y_1|\leq {\rm diam}\left\{y_1,..,y_k\right\}=1
\]
so another application of (\ref{e:exactfive}) and (\ref{e:exactseven})
yields $|\Phi_0(y_i)|\leq 2$. Hence,
$\Phi_0(x)\notin B_i^{+}$,
see (\ref{e:exacteight}). Consequently, (\ref{e:exacttwelve}) yields
\beq
(\Psi o\Phi_0)(x)=\Phi_0(x),\, |x|\geq 5.
\label{e:exactfifteen}
\eeq
From (\ref{e:exactthirteen}), we obtain that
\beq
\Psi o\Phi:\mathbb R^D\to \mathbb R^D
\label{e:exactsixteen}
\eeq
is an $\varepsilon$-distorted diffeomorphism since $\Phi_0$ is an Euclidean motion.
Next, applying Lemma~\ref{l:lemmatwo},
with $r_1=10$ and $r_2=\lambda^{-1/2}$,
we obtain an $\varepsilon$-distorted diffeomorphism
$\Psi:\mathbb R^D\to \mathbb R^D$ such that
\beq
\Psi_1(x)=\Phi_0(x),\, |x|\leq 10
\label{e:exactseventeen}
\eeq
and
\beq
\Psi_1(x)=x,\, |x|\geq \lambda^{-1/2}.
\label{e:exacteighteen}
\eeq

Note that Lemma~\ref{l:lemmatwo} applies, thanks to (\ref{e:exactseven}) and because we may assume
$\frac{\lambda^{-1/2}}{10}>\eta^{-1}$
with $\eta$ as in Lemma~\ref{l:lemmatwo}, thanks to our
``small $\lambda$ hypothesis''.

We now define
\begin{eqnarray}
\tilde{\Psi}(x):=
\left\{
\begin{array}{ll}
(\Psi o\Phi_0)(x), & |x|\leq 10 \\
\Psi_1(x), & |x|\geq 5.
\end{array}
\right\}.
\label{e:exactnineteen}
\end{eqnarray}

In the overlap region $5\leq |x|\leq 10$,
(\ref{e:exactfifteen}) and (\ref{e:exactseventeen}) show that
$(\Psi o\Phi_0)(x)=\Phi_0(x)=\Psi_1(x)$
so (\ref{e:exactnineteen}) makes sense.

We now check that $\tilde{\Psi}:\mathbb R^D\to \mathbb R^D$
is one to one and onto. To do so, we introduce the sphere
$S:=\left\{x:\, |x|=7\right\}\subset \mathbb R^D$ and partition $\mathbb R^D$ into $S$, inside($S$) and outside($S$).
Since $\Psi_1:\mathbb R^D\to \mathbb R^D$ is one to one and onto, (\ref{e:exactseventeen}) shows that the map
\beq
\Psi_1: {\rm outside}\, (S)\to {\rm outside}\, (\Phi_0(S))
\label{e:exacttwenty}
\eeq
is one to one and onto. Also, since $\Psi o\Phi_0:\mathbb R^D\to \mathbb R^D$ is one to one and onto, (\ref{e:exactfifteen}) shows that the map
\beq
\Psi o \Phi_0: {\rm inside}\,
(S)\to {\rm inside}\, (\Phi_0(S))
\label{e:exacttwentyone}
\eeq
is one to one and onto. In addition,
(\ref{e:exactfifteen}) shows that the map
\beq
\Psi o\Phi_0: (S)\to (\Phi_0(S))
\label{e:exacttwentytwo}
\eeq
is one to one and onto. Comparing (\ref{e:exactnineteen}) with (\ref{e:exacttwenty}), (\ref{e:exacttwentyone}) and
(\ref{e:exacttwentytwo}), we see that
$\tilde{\Psi}:\mathbb R^D\to \mathbb R^D$
is one to one and onto. Now since, also $\Psi o\Phi_0$ and $\Psi_1$ are
$\varepsilon$-distorted, it follows at once from (\ref{e:exactnineteen})
that $\tilde{\Psi}$ is smooth and
\[
(1+\varepsilon)^{-1}I\leq
(\tilde{\Psi'}(x))^{T}\tilde{\Psi'}(x)
\leq (1+\varepsilon)I,\, x\in \mathbb R^D.
\]
Thus,
\beq
\tilde{\Psi}:\mathbb R^{D}\to \mathbb R^{D}
\label{e:exacttwentythree}
\eeq
is an $\varepsilon$-distorted diffeomorphism.
From (\ref{e:exacteighteen}), (\ref{e:exactnineteen}),
we see that $\tilde{\Psi}(x)=x$
for $|x|\geq \lambda^{-1/2}$.
From ({\ref{e:exactforteen}), (\ref{e:exactnineteen}),
we have $\tilde{\Psi}(y_i)=z_i$ for each $i$, since, as we recall,
\[
|y_i|=|y_i-y_1|\leq {\rm diam}\left\{y_1,...,y_k\right\}=1.
\]
Thus, $\tilde{\Psi}$ satisfies
all the assertions in the statement of the
Theorem and the proof of the Theorem is complete.
$\Box$

\section{Proof of Theorem~\ref{t:lemmafive} and Theorem~\ref{t:lemmafiveprime}.}
\setcounter{equation}{0}

In this section, we prove Theorem~\ref{t:lemmafive} and Theorem~\ref{t:lemmafiveprime}. Theorem~\ref{t:lemmafive} follows immediately from the following result which is of independent interest.
Theorem~\ref{t:lemmafiveprime} follows exactly as in the proof of  Theorem~\ref{t:lemmafive} using Theorem~\ref{t:lemmafoura} in place of Theorem~\ref{t:lemmafour}.

\begin{thm}
Given $\varepsilon>0$, there exist $\eta,\delta>0$ depending on $\varepsilon$ such that the following holds. Let $E:=\left\{y_1,...,y_k\right\}$
and $E':=\left\{z_1,...,z_k\right\}$ be $k\geq 1$ distinct points of $\mathbb R^D$ with $1\leq k\leq D$ and $y_1=z_1$.
Suppose
\beq
(1+\delta)^{-1}\leq\frac{|z_i-z_j|}{|y_i-y_j|}\leq (1+\delta),\, i\neq j.
\label{e:lemmasixone}
\eeq
Then, there exists an $\varepsilon$-distorted diffeomorphism $\Phi:\mathbb R^D\to \mathbb R^D$ such that
\beq
\Phi(y_i)=z_i
\label{e:lemmasixtwo}
\eeq
for each $i$ and
\beq
\Phi(x)=x
\label{e:lemmasixthree}
\eeq
for
\[|x-y_1|\geq \eta^{-1}{\rm diam}\left\{y_1,...,y_k\right\}.\]
\label{l:lemmasix}
\end{thm}
\medskip

\noindent
{\bf Proof} \ We use induction on $k$. The case $k=1$ is trivial, we can just take
$\Phi$ to be the identity map. For the induction step, we fix $k\geq 2$ and suppose we already know the Theorem when $k$ is replaced by $k'<k$. We will prove the Theorem for the given $k$.
Let $\varepsilon>0$ be given. We pick small positive numbers $\delta'$, $\lambda$, $\delta$ as follows.
\beq
\delta'<\alpha(\varepsilon,D).
\label{e:lemmasixfour}
\eeq
\beq
\lambda<\alpha'(\delta',\varepsilon,D).
\label{e:lemmasixfive}
\eeq
\beq
\delta<\alpha''(\lambda,\delta',\varepsilon,D).
\label{e:lemmasixsix}
\eeq
Now let $y_1,...,y_k,z_1,...,z_k\in \mathbb R^D$ satisfy (\ref{e:lemmasixone}). We must produce an $\varepsilon$-distorted diffeomorphism $\Phi:\mathbb R^D\to\mathbb R^D$ satisfying
(\ref{e:lemmasixtwo}) and (\ref{e:lemmasixthree}) for some $\eta$ depending only on
$\delta,\lambda,\delta',\varepsilon, D$. That will complete the proof of the Theorem.

We apply Lemma~\ref{l:lemmathree} to $E=\left\{y_1,...,y_k\right\}$ with $\lambda$ in place of
$\eta$. Thus, we obtain an integer $l$ and a partition of $E$ into subsets
$E_1,E_2,...,E_{\mu_{\rm max}}$ with the following properties:
\beq
10\leq l\leq 100+\binom{k}{2}.
\label{e:lemmasixseven}
\eeq
\beq
{\rm diam}(E_{\mu})\leq \lambda^l{\rm diam}(E)
\label{e:lemmasixeight}
\eeq
for each $\mu$.
\beq
{\rm dist}(E_{\mu},E_{\mu'})\geq \lambda^{l-1}{\rm diam}(E)
\label{e:lemmasixnine}
\eeq
for $\mu\neq \mu'$. Note that
\beq
{\rm card}(E_{\mu})<{\rm card}(E)=k
\label{e:lemmasixten}
\eeq
for each $\mu$ thanks to (\ref{e:lemmasixeight}). For each $\mu$, let
\beq
I_{\mu}:=\left\{i:\, y_i\in E_{\mu}\right\}.
\label{e:lemmasixeleven}
\eeq

For each $\mu$, we pick a ``representative'' $i_{\mu}\in I_{\mu}$. The
$I_1,...,I_{\mu_{\rm max}}$ form a partition of $\left\{1,...,k\right\}$. Without loss of generality, we may suppose
\begin{equation}\label{eq lemmasixtwelve}
i_1=1.
\end{equation}
Define
\begin{eqnarray}\label{eq lemmasixthirteen}
&& I_{\rm rep}:=\left\{i_{\mu}:\mu=1,...,\mu_{\rm max}\right\} \\
&& E_{\rm rep}:=\left\{y_{i_{\mu}}:\, \mu=1,...,\mu_{\rm max}\right\}.
\end{eqnarray}
From (\ref{e:lemmasixeight}), (\ref{e:lemmasixnine}), we obtain
\[
(1-2\lambda^l){\rm diam}(E)\leq {\rm diam}(E_{\rm rep})\leq {\rm diam}(E),
\]
and
\[
|x'-x''|\geq \lambda^{l-1}{\rm diam}(E)
\]
for $x,x'\in S_{\rm rep}$, $x'\neq x''$.
Hence,
\beq
(1/2){\rm diam}(E)\leq {\rm diam}(E_{\rm rep})\leq {\rm diam}(E)
\label{e:lemmasixforteen}
\eeq
and
\beq
|x'-x''|\geq \lambda^m{\rm diam}(E_{\rm rep})
\label{e:lemmasixfifteen}
\eeq
for $x',x''\in E_{\rm rep},\, x'\neq x''$ where
\beq
m=100+\binom{D}{2}.
\label{e:lemmasixsixteen}
\eeq
See (\ref{e:lemmasixseven}) and recall that $k\leq D$.
We now apply Theorem~\ref{t:theorem3exact}
to the points $y_i,\,i\in I_{\rm rep}$, $z_i,\,i\in I_{\rm rep}$ with $\varepsilon$ in Theorem~\ref{t:theorem3exact} replaced by our present $\delta'$. The hypothesis of Theorem~\ref{t:theorem3exact} holds, thanks to the smallness assumptions (\ref{e:lemmasixfive}) and
(\ref{e:lemmasixsix}). See also
(\ref{e:lemmasixsixteen}), together with our present hypothesis (\ref{e:lemmasixone}).
Note also that $1\in I_{\rm rep}$ and $y_1=z_1$.
Thus we obtain a $\delta'$-distorted diffeomorphism $\Phi_0:\mathbb R^D\to \mathbb R^D$ such that
\beq
\Phi_0(y_i)=z_i,\, i\in I_{\rm rep}
\label{e:lemmasixseventeen}
\eeq
and
\beq
\Phi_0(x)=x \ \mbox{for} \ |x-y_1|\geq \lambda^{-1/2}{\rm diam}\left\{y_1,...,y_k\right\}.
\label{e:lemmasixeighteen}
\end{equation}
Define
\beq
y_i'=\Phi_0(y_i),\, i=1,...,k.
\label{e:lemmasixnineteen}
\eeq
Thus,
\beq
y_{i_{\mu}}'=z_{i_{\mu}}
\label{e:lemmasixtwenty}
\eeq
for each $\mu$ and
\beq
(1+C\delta')^{-1}\leq \frac{|z_i-z_j|}{|y_i'-y_j'|}\leq (1+C\delta'),\, i\neq j
\label{e:lemmasixtwentyone}
\eeq
thanks to (\ref{e:lemmasixone}), (\ref{e:lemmasixsix}), (\ref{e:lemmasixnineteen})
and the fact that $\Phi_0:\mathbb R^D\to \mathbb R^D$ is a $\delta'$-distorted diffeomorphism.
Now fix $\mu (1\leq \mu\leq \mu_{{\rm max}})$.
We now apply our inductive hypothesis with $k'<k$ to the points $y_i',z_i,\, i\in I_{\mu}$. (Note that the inductive
hypothesis applies, thanks to (\ref{e:lemmasixten})
Thus, there exists
\beq
\eta_{\rm indhyp}(D, \varepsilon)>0,\, \delta_{\rm indhyp}(D, \varepsilon)>0
\label{e:lemmasixtwentytwo}
\eeq
such that the following holds: Suppose
\beq
(1+\delta_{\rm indhyp})^{-1}|y_i'-y_j'|\leq |z_i-z_j|\leq |y_i'-y_j'|(1+\delta_{\rm indhyp}),\, i,j\in I_{\mu}
\label{e:lemmasixtwentythree}
\eeq
and
\beq
y_{i_{\mu}}'=z_{i_{\mu}}.
\label{e:lemmasixtwentyfour}
\eeq
Then there exists a $\varepsilon$ distorted diffeomorphism $\Psi_{\mu}:\mathbb R^D\to \mathbb R^D$ such that
\beq
\Psi_{\mu}(y_i')=z_i,\, i\in I_{\mu}
\label{e:lemmasixtwentyfive}
\eeq
and
\beq
\Psi_{\mu}(x)=x,\ \mbox{for} \ |x-y_{i_{\mu}}'|\geq \eta_{\rm indhyp}^{-1}{\rm diam}(S_{\mu}).
\label{e:lemmasixtwentysix}
\eeq

We may suppose $C\delta'<\delta_{\rm indhyp}$
with $C$ as in (\ref{e:lemmasixtwentyone}),
thanks to (\ref{e:lemmasixtwentytwo}) and our smallness assumption
(\ref{e:lemmasixfour}). 
Similarly, we may suppose that $\eta_{\rm indhyp}^{-1}<1/2\lambda^{-1/2}$,
thanks to (\ref{e:lemmasixtwentytwo}) and our smallness assumption
(\ref{e:lemmasixfive}). Thus (\ref{e:lemmasixtwentythree}) and (\ref{e:lemmasixtwentyfour}) hold, by virtue of (\ref{e:lemmasixtwentyone}) and (\ref{e:lemmasixtwenty}). Hence, for each $\mu$, we obtain an $\varepsilon$-distorted diffeomorphism $\Psi_{\mu}:\mathbb R^{D}\to \mathbb R^{D}$, satisfying
(\ref{e:lemmasixtwentyfive}) and (\ref{e:lemmasixtwentysix}). In particular,
(\ref{e:lemmasixtwentysix}) yields
\beq
\Psi_{\mu}(x)=x,\ \mbox{for} \ |x-y_{i_{\mu}}'|\geq 1/2\lambda^{-1/2}{\rm diam}(S_{\mu}).
\label{e:lemmasixtwentyseven}
\eeq
Taking
\beq
B_{\mu}=B(y_{i_{\mu}}^\prime, 1/2\lambda^{-1/2}{\rm diam}(S_{\mu}) ),
\label{e:lemmasixtwentyeight}
\eeq
we see from (\ref{e:lemmasixtwentyseven}), that
\beq
\Psi_{\mu}:B_{\mu}\to B_{\mu}
\label{e:lemmasixtwentynine}
\eeq
is one to one and onto since $\Psi_{\mu}$ is one to one and onto.
Next, we note that the balls $B_{\mu}$ are pairwise disjoint.* (Note that the closed ball $B_{\mu}$
is a single point if $S_{\mu}$ is a single point.) This follows from (\ref{e:lemmasixeight}), (\ref {e:lemmasixnine}) and the definition (\ref{e:lemmasixtwentyeight}).
We may therefore define a map $\Psi:\mathbb R^D\to \mathbb R^D$ by setting
\beq
\Psi(x):=\left\{
\begin{array}{ll}
\Psi_{\mu}(x), & x\in B_{\mu}\, , \ {\rm any} \ \mu \\
x, & x\notin \cup_{\mu} B_{\mu}
\end{array}
\right\}.
\label{e:lemmasixthirty}
\eeq
Thanks to (\ref{e:lemmasixtwentynine}), we see that $\Psi$ maps $\mathbb R^D$ to $\mathbb R^D$
one to one and onto. Moreover, since each $\Psi_{\mu}$ is an $\varepsilon$-distorted diffeomorphism
satisfying (\ref{e:lemmasixtwentyseven}), we see that $\Psi$ is smooth on $\mathbb R^D$ and that
\[
(1+\varepsilon)^{-1}I\leq (\Psi'(x))^{T}\Psi'(x)\leq (1+\varepsilon)I,\, x\in \mathbb R^D.
\]
Thus,
\beq
\Psi:\mathbb R^D\to \mathbb R^D
\label{e:lemmasixthirtyone}
\eeq
is an
$\varepsilon$-distorted diffeomorphism. From (\ref{e:lemmasixtwentyfive}) and
(\ref{e:lemmasixthirty}), we see that
\beq
\Psi(y_i')=z_i,\, i=1,...,k.
\label{e:lemmasixthirtytwo}
\eeq
Let us define
\beq
\Phi=\Psi o\Phi_0.
\label{e:lemmasixthirtythree}
\eeq
Thus
\beq
\Phi \,
\mbox{is a} \ C \varepsilon \ \mbox{-distorted diffeomorphism of} \
\mathbb R^D\to \mathbb R^D
\label{e:lemmasixthirtyfour}
\eeq
since $\Psi,\Phi_0:\mathbb R^D\to \mathbb R^D$ are
$\varepsilon$ distorted diffeomorphisms.
Also
\beq
\Phi(y_i)=z_i,\, i=1,...,k
\label{e:lemmasixthirtyfive}
\eeq
as we see from (\ref{e:lemmasixnineteen}) and (\ref{e:lemmasixthirtytwo}). Now suppose that
\[
|x-y_1|\geq \lambda^{-1}{\rm diam}\left\{y_1,...,y_k\right\}.
\]
Since $\Phi_0:\mathbb R^D\to \mathbb R^D$ is a $\varepsilon$-distorted diffeomorphism, we have
\beq
\left|\Phi_0(x)-y_1'\right|\geq (1+\varepsilon)^{-1}\lambda^{-1}{\rm diam}\left\{y_1,...,y_k\right\}
\label{e:lemmasixthirtysix}
\eeq
and
\[
{\rm diam}\left\{y_1',...,y_k'\right\}\leq (1+\varepsilon){\rm diam}\left\{y_1,...,y_k\right\}.
\]
See (\ref{e:lemmasixnineteen}).

Hence for each $\mu$,
\begin{eqnarray*}
&& \left|\Phi_0(x)-y_{i_{\mu}}'\right|\geq \left[(1+\varepsilon)^{-1}\lambda^{-1}-(1+\varepsilon)
\right]{\rm diam}\left\{y_1,...,y_k\right\} \\
&& >1/2\lambda^{-1/2}{\rm diam}(S_{\mu}).
\end{eqnarray*}
Thus, $\Phi_0(x)\notin \cup_{\mu}B_{\mu}$, see (\ref{e:lemmasixtwentyeight}), and therefore
$\Psi o\Phi_{0}(x)=\Phi_{0}(x)$, see (\ref{e:lemmasixthirty}). Thus,
\beq
\Phi(x)=\Phi_0(x).
\label{e:lemmasixthirtyseven}
\eeq
From (\ref{e:lemmasixeighteen}) and
(\ref{e:lemmasixthirtyseven}), we see that $\Phi(x)=x$.
Thus, we have shown that
\[
|x-y_1|\geq \lambda^{-1}{\rm diam}\left\{y_1,...,y_k\right\}
\]
implies $\Phi(x)=x$. That is,
(\ref{e:lemmasixthree}) holds with $\eta=\lambda$. Since
also (\ref{e:lemmasixthirtyfour}) and (\ref{e:lemmasixthirtyfive}) hold we have carried out
our inductive step completely and hence, the proof of the Theorem $\Box$.

In this last section, we address the restriction $k\leq D$ in Theorem~\ref{t:lemmafive}. 

\section{A Counterexample.}
\setcounter{equation}{0}

In this final section, we show that in general one cannot expect to have
Theorem~\ref{t:lemmafive} and Theorem~\ref{t:lemmafiveprime} for $k>D$. This is shown using the following counterexample.
\medskip

\noindent
{\bf Counterexample:} \ Fix $2D+1$ points as follows. Let $\delta>0$ be a small
positive number depending on $D$. Let $y_1,...,y_{D+1}\in \mathbb R^D$
be the vertices of a regular simplex, all lying on the
sphere of radius $\delta$ about the origin.
Then define $y_{D+2} \cdots y_{2D+1}\in \mathbb R^D$ such that
$y_{D+1},...,y_{2D+1}$ are the vertices of a regular simplex, all lying in a sphere of radius 1, centered at some point $w_1\in \mathbb R^D$. Next, we define a map
\[
\phi:\left\{y_1,...,y_{2D+1}\right\}\to \left\{y_1,...,y_{2D+1}\right\}
\]
as follows. We take $\phi|_{\left\{y_1,...,y_{D+1}\right\}}$
to be an odd permutation that fixes
$y_{D+1}$, and take $\phi|_{\left\{y_{D+1},...,y_{2D+1}\right\}}$ to be the identity.
The map $\phi$ distorts distances by at most a factor $1+C\delta$. Here, we can take $\delta$
arbitrarily small. On the other hand, for small enough $\varepsilon_0>0$ depending only on $D$, we will show that $\phi$ cannot be extended to a map $\Phi:\mathbb R^D\to \mathbb R^D$ satisfying
\[
(1+\varepsilon_0)^{-1}|x-x'|\leq|\Phi(x)-\Phi(x')|\leq |x-x'|(1+\varepsilon_0),\, x,x'\in \mathbb R^D.
\]
In fact, suppose that such a $\Phi$ exists.
Then $\Phi$ is continuous. Note that there exists
$T\in O(D)$ with ${\rm det}T=-1$ such that $\phi(y_i)=Ty_i$ for $i=1,...,D+1$. Let
$S_t$ be the sphere of radius $r_t:=\delta\cdot(1-t)+1\cdot t$ centered at
$t \cdot w_1$ for $t\in [0,1]$ and let
$S_t'$ be the sphere of radius $r_t$ centered at
$\Phi(t \cdot w_1)$. Also, let $Sh_t$ be the spherical shell
\[
\left\{x\in \mathbb R^D:\, r_t\cdot(1+\varepsilon_0)^{-1}\leq |x-
\Phi(t \cdot w_1)|\leq r_t\cdot(1+\varepsilon_0)\right\}
\]
and let $\pi_t:Sh_t\to S_t'$ be the projection
defined by
\[
\pi_t(x)-\Phi(t \cdot w_1)=
\frac{x-\Phi(t \cdot w_1)}{|x-\Phi(t \cdot w_1)|} \cdot r_t.
\]
Since $\Phi$ agrees with $\phi$, we know that
\beq
\left|\Phi(x)-Tx\right|\leq C\varepsilon_0\delta,\, |x|=\delta.
\label{e:counterone}
\eeq
Since $\Phi$ agrees with $\phi$, we know that
\beq
\left|\Phi(x)-x\right|\leq C\varepsilon_0,\, |x-w_1|=1.
\label{e:countertwo}
\eeq
Our assumption that $\Phi$ is an approximate isometry shows that
\[
\Phi: S_t\to Sh_t,\, 0\leq t\leq 1
\]
and
\beq
(\pi_t)o(\Phi):S_t\to S_t',\, 0\leq t\leq 1.
\label{e:counterthree}
\eeq
We can therefore define a one-parameter family of maps $\Psi_t,\, t\in [0,1]$ from the unit sphere
to itself by setting
\begin{eqnarray*}
&& \Psi_t(w)=\frac{(\pi_t o\Phi)(tw_1+r_tw)-\Phi(tw_1)}
{\left|(\pi_t o\Phi)(tw_1+r_tw)-\Phi(tw_1)\right|}= \\
&& \frac{(\pi_t o\Phi)(tw_1+r_tw)-\Phi(tw_1)}{r_t}
\end{eqnarray*}
from the unit sphere to itself. Then $\Psi_t$ is a continuous family of continuous maps from the unit sphere to itself. From (\ref{e:counterone}), we see that $\Psi_0$ is a small perturbation of the map $T:S^{D-1}\to S^{D-1}$ which has degree -1. From (\ref{e:countertwo}), we see that
$\Psi_1$ is a small perturbation of the identity.
Consequently, the following must hold:
\begin{itemize}
\item Degree $\Psi_t$ is independent of $t\in [0,1]$.
\item Degree $\Psi_0=-1$.
\item Degree $\Psi_1=+1$.
\end{itemize}
Thus, we have arrived at a contradiction, proving that there does not exist an extension $\Phi$ as above. $\Box$
\medskip

The counter example above could be motivated by the fact that the mapping swapping the real numbers $\eta$ and $-\eta$ and fixing the number $1$ cannot be extended to a continuous bijection of the line.
\medskip

The case of $k>D$ absent in Theorem~\ref{t:lemmafive} and Theorem~\ref{t:lemmafiveprime} is investigated by us in the paper \cite{FD} where it is shown that we may allow for this case if roughly we require that on any $D+1$ of the $k$ points which form a relatively voluminous simplex, the extension $\Phi$ is orientation preserving.
\medskip

\noindent
{\bf Acknowledgment:}\, Support from the Department of Mathematics at Princeton, the National Science Foundation, The American Mathematical Society, The Department of Defense and the University of the Witwatersrand are gratefully acknowledged.

\end{document}